\documentclass{article}
\usepackage{amsmath, amsfonts, amsthm, amssymb,fullpage}
\usepackage{comment}
\usepackage{hyperref}
\usepackage{url}
\usepackage{color}
\definecolor{mylinkcolor}{rgb}{0.5,0.0,0.0}
\definecolor{myurlcolor}{rgb}{0.0,0.0,0.75}
\hypersetup{colorlinks=true,urlcolor=myurlcolor,citecolor=myurlcolor,linkcolor=mylinkcolor,linktoc=page,breaklinks=true}

\def\Sage{{\tt Sage}}
\def\Magma{{\tt Magma}}
\def\Pari{{\tt Pari}}
\def\Z{{\mathbb Z}}
\def\Q{{\mathbb Q}}

\def\F{{\mathbb F}}
\def\Fp{{\mathbb F}_p}
\def\OO{{\mathcal O}}
\def\a{{\mathfrak a}}
\def\b{{\mathfrak b}}
\def\idc{{\mathfrak c}}
\def\p{{\mathfrak p}}
\def\q{{\mathfrak q}}
\def\n{{\mathfrak n}}

\DeclareMathOperator{\disc}{disc}
\DeclareMathOperator{\ord}{ord}

\title{Sorting and labelling integral ideals in a number field}
\author{John Cremona, Aurel Page, and Andrew V. Sutherland}

\begin{document}
\maketitle

We define a scheme for labelling and ordering integral ideals of
number fields, including prime ideals as a special case.  The order we
define depends only on the choice of a monic irreducible integral
defining polynomial for each field, and we start by defining for each
field its unique \emph{reduced} defining polynomial in
Section~\ref{sec:polredabs}, after Belabas.  In
Section~\ref{sec:primes} we define a total order on the set of prime
ideals of~$K$ and then extend this in Section~\ref{sec:ideals} to a
total order on the set of all nonzero integral ideals of~$K$.  For all
integral ideals~$\a$, including prime ideals, this order allows us to
give a unique label of the form $N.i$, where~$N=N(\a)$ and $i\ge1$ is
the index of~$\a$ in the ordered list of all ideals of norm~$N$.  Our
ideal labelling scheme has several nice properties: for a given norm,
prime ideals always appear first, and given the factorisation of~$N$,
the bijection between ideals of norm~$N$ and labels is computable in
polynomial time.

Our motivation for this is to have a well-defined and concise way to
sort and label ideals for use in databases such as the
LMFDB~\cite{lmfdb}.  For example, consider the finite set of Bianchi
modular newforms defined over the imaginary quadratic field~$K$ having
fixed level~$\n$, weight~$2$ and trivial character.  In order to label
these, we sort them according to the vector of their Fourier
coefficients~$a_{\p}$, indexed by prime ideals~$\p$ of~$K$, and for
this to be well-defined we need to specify an order on the primes
themselves.

We have implemented algorithms which realise this scheme, in \Sage,
\Magma\ and~\Pari. Each of the three authors wrote one version of the code
independently of the others, and we checked on a series of examples that the
outputs agree. The code is available
at~\url{https://github.com/JohnCremona/Sorting/tree/master/code}.

\section{Defining the number field}
\label{sec:polredabs}
Every number field has the form $K=\Q(\alpha)$ where $\alpha$ is an
algebraic integer.  The minimal polynomial of~$\alpha$ is an
irreducible monic polynomial in~$\Z[X]$, called the \emph{defining
  polynomial} of~$K$.  While every field has infinitely many such
defining polynomials, for our purposes, we need to make a choice once
and for all of one defining polynomial for each field, since the way
in which we order field elements and ideals will depend on this
choice.

The \emph{reduced defining polynomial} of the number field~$K$ is
defined to be the unique monic integral polynomial determined by the
following process. Recall that the $T_2$-norm of a polynomial is the
sum of the squares of the modulus of the roots, and for a
polynomial~$P = x^n + a_1x + \dots + a_n$, let~$S(P) =
(|a_1|,a_1,\dots,|a_n|,a_n)$; there are only finitely many monic integral
polynomials of a given degree and bounded~$T_2$ norm.
\begin{enumerate}
  \item Let~$L_0$ be the finite list of monic integral defining polynomials
    of~$K$ which are minimal with respect to the~$T_2$ norm.
  \item Let~$L_1$ be the sublist of~$L_0$ of polynomials whose
    discriminant has minimal absolute value.
  \item Order the polynomials in~$L_1$ by lexicographic order of the
    vectors~$S(P)$.

    Then the reduced defining polynomial of~$K$ is the smallest
    polynomial in~$L_1$ with respect to this order.
\end{enumerate}

Note that in the second step, we restrict to the list~$L_0$ of
polynomials defined in the first step; hence the defining polynomial
for the field with smallest discriminant may not belong to the list.

The third step distinguishes, for example, between $x^2-x+1$ and
$x^2+x+1$ as defining polynomials for $\Q(\sqrt{-3})$, choosing the
first, and hence the generator $(1+\sqrt{-3})/2$ (a 6th root of unity)
rather than $(-1+\sqrt{-3})/2$ (a 3rd root of unity).

This definition is due to Karim Belabas, and is implemented in the
\Pari\ function {\tt polredabs()}, which is guaranteed to always
return the same polynomial, given any defining polynomial for a number
field: see \cite{polredabs}.  Note that computing this polynomial can
be very time-consuming when the degree or discriminant of the field is
large: in particular it requires computation of the ring of integers
of~$K$, and hence the factorization of the discriminant of an initial
(monic integral irreducible) defining polynomial for~$K$; in addition the number
of potential defining polynomials to consider could be exponential in the
degree.

We do not assume that $\Z[\alpha]$ is the full ring of integers
$\OO_K$.  Indeed, for most number fields $K$ this condition is
not satisfied by any choice of $\alpha$ (in general the ring $\OO_K$
will not be monogenic over $\Z$).

We emphasise that both the defining polynomial $g(X)$ and the
generator $\alpha$ of~$K$ are fixed from now on.  In other words we
are fixing the structure of $K=\Q[X]/(g(X))$ not just as a field, but
as a $\Q[X]$-algebra, where the structure map $\Q[X]\to K$ has kernel
$(g(X))$ and $\alpha$ is the image of $X$ in $K$.  In what follows,
the order we define on field elements and ideals depends on the choice
of defining polynomial, so is canonical when the reduced defining
polynomial is used.

\section{Sorting prime ideals}\label{sec:primes}
Fix a number field $K=\Q(\alpha)=\Q[X]/(g(X))$ with defining
polynomial $g\in\Z[X]$ and generator $\alpha$ as above, and let
$\OO_K$ denote its ring of integers.  We do not assume that $K/\Q$ is
Galois or that $\OO_K=\Z[\alpha]$.  By a \emph{prime} of $K$ (or of
$\OO_K$) we mean a nonzero prime ideal in the ring $\OO_K$; for primes
of $\Q$ we may also use $p$ to denote the prime ideal $(p)$.

We wish to define a total order on the set of prime ideals~$\p$
of $\OO_K$.  Let $(p)=\p\cap \Z$ be the prime of $\Z$ lying below $\p$,
let $e=e(\p/p)$ be the ramification index (the multiplicity of $\p$ in the prime factorization of $(p)$ in the Dedekind domain $\OO_K$), and let $f=f(\p/p)$ the residue field degree $[\OO_K/\p:\F_p]$.
The (absolute) norm of $\p$ is then $N(\p)=p^f$.


We order the primes $\p$ first by norm and then by the ramification
index.  It remains only to order the primes~$\p$ of $K$ that lie above
the same rational prime $p$ and have the same norm (hence the same
residue degree) and ramification index.  Such primes will necessarily
have the same value of $p^{ef}$; the converse need not hold (primes
with the same value of $p^{ef}$ need not have the same norm $p^{f}$),
but this will not concern us since we order by norm first.

\subsection{The case of unramified primes}

Our general method is simplest when $p$ is unramified, not just in $K$
but in the order $\Z[\alpha]$; this holds precisely when $p\nmid
[\OO_K:\Z[\alpha]]\disc(K)$, or equivalently, $p\nmid\disc(g)$, which
applies to all but finitely many $p$.  In this case $g(X)$ factors
modulo~$p$ into distinct irreducible elements of $\Fp[X]$, and we may
write
\[
    g(X) \equiv \prod_i h_i(X) \pmod{p},
\]
with the factors $h_i(X)\in\Z[X]$ monic and such that their reductions
modulo~$p$ are distinct and irreducible.  We may assume that the
coefficients of each $h_i(X)$ are reduced modulo~$p$ to integers in the interval $[0,p-1]$, and we order these factors first by degree and then lexicographically by their coefficient vectors in $[0,p-1]^{f+1}$.

The Dedekind-Kummer theorem allows us to associate to each factor
$h_i(X)$ the prime ideal $\p_i=(p,h_i(\alpha))$ above $p$, with
residue degree $f(\p_i/p)=\deg h_i$ and norm $N(\p_i)=p^{\deg h_i}$.
We then have $p\OO_K=\prod_i \p_i$ and our ordering of the $h_i$
induces an ordering of the primes~$\p_i$ above $p$.  If we have any
representation of $\p_i$ in the form $\p_i=(p,\beta)$ where
$\beta=b(\alpha)$ for some $b(X)\in\Z[X]$, then we can recover
$h_i(X)$ since $h_i(X)=\gcd(g(X),b(X))$ (with the $\gcd$ computed in
$\F_p[X]$).  This follows from the observation that
$\ord_{\p_i}(\beta)>0$ while $\ord_{\p_j}(\beta)=0$ for all primes
$\p_j$ above~$p$ other than~$\p_i$.  This makes implementation of the
sorting function on the primes above~$p$ very simple in this case.

\subsection{The general case, including ramified primes}

There is a bijection between the distinct primes $\p$ above $p$ and
the irreducible factors $h(X)$ of $g(X)$ in $\Q_p[X]$ in which $\deg h
= e(\p/p)f(\p/p)$ (see \cite[Theorem 3.8 (d)]{Janusz}).  Write
\[
  g(X) = h_1(X)h_2(X)\dots h_r(X)
\]
with the $h_i(X)$ monic and irreducible in $\Z_p[X]$. Since we sort
prime ideals by norm and ramification index, it is enough to describe
an order on polynomials of the same degree; in fact, we only need to
order the set of polynomials $h_i(X)$ that have the same degree~$ef$,
and also correspond to primes $\p_i$ with both the same ramification
index~$e$ and the same residue degree~$f$, which may be a smaller set.

We sort polynomials in~$\Z_p[X]$ of degree~$d$ as follows. Such a polynomial~$P$
can be written uniquely in the form
\[
  (a_{0,0}+a_{0,1}p+a_{0,2}p^2+\dots) + (a_{1,0}+a_{1,1}p+\dots)X
  + \dots + (a_{d,0}+a_{d,1}p+\dots)X^d,
\]
with all~$a_{i,j}\in\Z\cap[0,p-1]$. We can then attach to~$P$ the infinite
vector
\[
(a_{0,0},a_{1,0},\dots,a_{d,0},a_{0,1},a_{1,1},\dots,a_{d,1},a_{0,2},\dots)\in [0,p-1]^\infty,
\]
and we sort polynomials according to the lexicographic ordering of such vectors.
Comparisons can be made using finite precision, provided that the polynomials are known to be
distinct, which is the case for the~$h_i$.


We now explain how to compute the bijection between the~$h_i$ and the primes
above~$p$. Let~$h\in\Z_p[X]$ be one of the~$h_i$, and let~$\p$ be a prime
above~$p$. The valuation~$v_{\p}$ on the number field~$K$ (a finite \'etale $\Q$-algebra) extends to the \'etale $\Q_p$-algebra ~$K\otimes_{\Q}\Q_p = \prod_i
\Q_p[X]/(h_i(X))$ obtained via base-change in the following way.
Let~$j$ be such that~$h_{j}$ corresponds
to~$\p$, and let~$v_j$ be the valuation of the $p$-adic field~$\Q_p[X]/(h_j(X))$.
Define the extension of~$v_\p$ to $K\otimes_\Q\Q_p$ to be the composition of the maps
\[
  K\otimes_{\Q}\Q_p \longrightarrow \Q_p[X]/(h_j(X)) \stackrel{v_j}{\longrightarrow} \Z\cup\{\infty\}.
\]
We then have~$v_{\p}(h(\alpha)) = \infty$ if and only
if~$\p$ corresponds to~$h$.
While this valuation cannot be computed using finite
approximations to~$h$, for all integers~$k\ge 1$ we have~$v_\p((p^k,h(\alpha)))
= \min(e(\p/p)k, v_\p(h(\alpha)))$, and this value equals~$e(\p/p)k$
if~$h$ corresponds to~$\p$ and is bounded above otherwise (independently of~$k$). The
valuation~$v_\p((p^k,h(\alpha)))$ can be computed using finite approximations
to~$h$, and if $k$ is such that the set $\{v_\p((p^k,h_i(\alpha)))\}$ has a unique maximum, this maximum occurs for the polynomial $h_i(X)$ corresponding to $\p$.
Thus we can compute the bijection by making $k$ sufficiently large (and any $k$ that yields a unique maximum works).

We remark that ``Round 4" of the $p$-adic polynomial factorization
algorithm implemented in \Pari\ \cite{Roblot} computes the prime
ideals corresponding to the various $p$-adic factors, but this part of
the output is not available via the standard interface.  The same is
true in \Sage, which uses the \Pari\ library for $p$-adic
factorization.  Nevertheless, the discussion above allows us to
recover the bijection without requiring direct access to this
implementation, which may not be the same in other computer algebra
systems in any case.

\subsection{Examples}
We use number fields in the LMFDB database (see~\cite{lmfdb}) to
illustrate the sorting of prime ideals.

Let $K=\Q(\alpha)$ be the non-Galois cubic field with LMFDB
label~\href{www.lmfdb.org/NumberField/3.1.503.1}{\texttt{3.1.503.1}}
of discriminant $-503$ with reduced defining polynomial
$g(x)=X^3-X^2+2X+8$.  The prime $2$ is unramified, but it divides
$\disc(g(X))=-2^2\cdot 503$ and is an essential divisor of
$[\OO_K:\Z[\alpha]]$ in the sense that $2\mid[\OO_K:\Z[\alpha']]$ for
every algebraic integer $\alpha'$ in $K$ (this example, due to
Dedekind, is the standard example of a non-monogenic field).

Let $p=2$ and $k=2$; the $2$-adic factors of $g(X)$ are: $h_1=X+O(2^2)$,
$h_2=X+2+O(2^2)$, $h_3=X+1+O(2^2)$.  So there are $3$
primes above $2$, each with residue degree $1$.
Calling these $\p_a = (2,\frac 12\alpha^2+\frac 12\alpha+3)$, $\p_b =
(2,\alpha+3)$, $\p_c = (2,\frac 12\alpha^2-\frac 12\alpha)$ in random order we find that
\begin{itemize}
  \item the $\p_a$-valuations of the $(2^2,h_i(\alpha))$ are $(1,2,0)$,
  \item the $\p_b$-valuations of the $(2^2,h_i(\alpha))$ are $(0,0,2)$,
  \item the $\p_c$-valuations of the $(2^2,h_i(\alpha))$ are $(2,1,0)$,
\end{itemize}
thus $\p_1=\p_c$, $\p_2=\p_a$ and $\p_3=\p_b$.

Let $p=503$.  There are two primes above~$p$, both of norm~$p$, only
one of which is ramified (with $e=2$). Our ordering puts the
unramified prime first and the ramified prime second, so we do not even need to
look at the $p$-adic factorization of $g(X)$ in this case.  In fact,
with $k=2$ we find factors $X+191929+O(503^2)$ and $X^2+61079X+87617+O(503^2)$.

For a larger example, let $K=\Q(\alpha)=\Q[X]/(g(X))$, where
\[
g(X)= X^{10} - 3X^9 - 35X^8 + 120X^7 + 242X^6 - 1080X^5 + 44X^4 + 2343X^3 -
1631X^2 + 111X + 79.
\]
This field has LMFDB
label~\href{http://beta.lmfdb.org/NumberField/10.10.24952891341003125.1}{\texttt{10.10.24952891341003125.1}} and discriminant~$5^5 41^8$, but we have
\[
\disc(g)=3^{12}\cdot 5^5\cdot 41^8\cdot 2141^2\cdot 26641^2.
\]
Let us consider the primes above $3$.
Over~$\Q_3$, the polynomial $g(X)$ splits as
a product of~$5$ polynomials of degree~$2$, that is, $g = h_1h_2h_3h_4h_5$ with
$h_i\in\Z_3[x]$. Modulo~$3$, we have
\begin{itemize}
  \item $h_1 = X^2 + 1 +O(3)$,
  \item $h_2 = X^2 + 2X + 2 + O(3)$,
  \item $h_3 = X^2 + X + 2 + O(3)$,
  \item $h_4 = X^2 + 2X + 2 + O(3)$,
  \item $h_5 = X^2 + X + 2 + O(3)$.
\end{itemize}
The prime $3$ is unramified, so we can already see that it splits into 5 primes of norm $9$.  These approximations are not sufficient to distinguish all the $h_i$,
but we get the following initial segments of the associated vectors (omitting
the maximal degree term since they are all monic):
\begin{itemize}
  \item $(1,0,\dots)$,
  \item $(2,2,\dots)$,
  \item $(2,1,\dots)$,
  \item $(2,2,\dots)$,
  \item $(2,1,\dots)$,
\end{itemize}
so we have~$h_1 < \{h_3,h_5\} < \{h_2,h_4\}$. Modulo~$9$, we get
\begin{itemize}
  \item $h_1 = X^2 + 3X + 1 + O(3^2)$,
  \item $h_2 = X^2 + 5X + 5 + O(3^2)$,
  \item $h_3 = X^2 + 7X + 2 + O(3^2)$,
  \item $h_4 = X^2 + 5X + 5 + O(3^2)$,
  \item $h_5 = X^2 + 4X + 5 + O(3^2)$.
\end{itemize}
This is still not enough to distinguish them, but
we get the refined initial segments
\begin{itemize}
  \item $(1,0,0,1,\dots)$,
  \item $(2,2,1,1,\dots)$,
  \item $(2,1,0,2,\dots)$,
  \item $(2,2,1,1,\dots)$,
  \item $(2,1,1,1,\dots)$,
\end{itemize}
so we have~$h_1<h_3<h_5<\{h_2,h_4\}$. Finally, modulo~$27$ we have
\begin{itemize}
  \item $h_1 = X^2 + 3X + 1 + O(3^3)$,
  \item $h_2 = X^2 + 5X + 5 + O(3^3)$,
  \item $h_3 = X^2 + 7X + 11 + O(3^3)$,
  \item $h_4 = X^2 + 23X + 23 + O(3^3)$,
  \item $h_5 = X^2 + 13X + 14 + O(3^3)$.
\end{itemize}
This is now enough to distinguish all the~$h_i$, and we obtain the initial segments
\begin{itemize}
  \item $(1,0,0,1,0,0\dots)$,
  \item $(2,2,1,1,0,0,\dots)$,
  \item $(2,1,0,2,1,0,\dots)$,
  \item $(2,2,1,1,2,2,\dots)$,
  \item $(2,1,1,1,1,1,\dots)$,
\end{itemize}
yielding the order~$h_1 < h_3 < h_5 < h_2 < h_4$.
As noted above, the prime~$3$ decomposes as a product of~$5$ prime ideals of norm~$9$,
which we denote~$\p_a, \dots, \p_e$ in arbitrary order. We know that we need
precision at least~$O(3^3)$ to distinguish the polynomials~$h_i$, so we compute
the valuations~$v_\p((3^3,h_i(\alpha)))$. We obtain:
\begin{itemize}
  \item The~$\p_a$-valuations of~$(3^3,h_i(\alpha))$ are~$(0,0,0,0,3)$,
  \item The~$\p_b$-valuations of~$(3^3,h_i(\alpha))$ are~$(0,3,0,2,0)$,
  \item The~$\p_c$-valuations of~$(3^3,h_i(\alpha))$ are~$(3,0,1,0,0)$,
  \item The~$\p_d$-valuations of~$(3^3,h_i(\alpha))$ are~$(0,2,0,3,0)$,
  \item The~$\p_e$-valuations of~$(3^3,h_i(\alpha))$ are~$(1,0,3,0,0)$,
\end{itemize}
thus we have the bijection~$\p_a\leftrightarrow h_5$, $\p_b\leftrightarrow h_2$,
$\p_c\leftrightarrow h_1$, $\p_d\leftrightarrow h_4$, $\p_e\leftrightarrow h_3$.
Since these ideals all have norm~$9$ and ramification index~$1$, our ordering of
the primes above~$3$ is~$\p_c < \p_e < \p_a < \p_b < \p_d$.

\section{Sorting all nonzero integral ideals}
\label{sec:ideals}
Fix a number field $K=\Q(\alpha)$ with defining polynomial
$g(X)\in\Z[X]$ and generator $\alpha$ as above.  We now define a total
order on the set of all nonzero integral ideals of~$K$.  We first
order ideals by norm, and then sort ideals of the same norm according
to the criterion which we specify in this section. This defines for
each ideal~$\a$ a unique label~$N.i$, where~$N$ is the norm of~$\a$
and~$i\ge 1$ is the index of~$\a$ in the ordered list of ideals of
norm~$N$.  This order has the following properties:
\begin{itemize}
  \item prime ideals are smaller than every non-prime ideal of the same norm,
    and are ordered in the same way as before;
  \item if~$\a,\b,\idc$ are integral ideals such that~$\a<\b$,
    then~$\a\idc<\b\idc$; so our ordering makes the set of nonzero
    integral ideals into a totally ordered monoid;
  \item if the factorisation of a positive integer~$N$ is known, then the
    bijection between ideals of norm~$N$ and their labels is computable in
    polynomial time (even though the number of such ideals might not be
    polynomial in~$\log N$).
\end{itemize}

\subsection{Ideals of prime power norm}\label{sec:primepowernorm}

We first define an order on the set of ideals of prime power
norm. Let~$\a$ be an ideal of norm~$p^n$, and let~$\p_1,\dots,\p_r$ be
the prime ideals of~$K$ above~$p$, ordered as in
Section~\ref{sec:primes}. We have~$\a = \p_1^{v_1}\dots \p_r^{v_r}$
for some integers~$v_i\ge 0$. We define the weight of this
factorization to be~$v_1+\dots+v_r$. We order ideals of norm a power
of $p$ by increasing norm first, then increasing weight, and finally
by reverse lexicographic order of the vector of
exponents~$(v_1,\dots,v_r)$. For this order, the prime ideals come
first (having weight~$1$), and are ordered in the same way as
previously.

Example: let~$K$ be a number field and~$p$ a prime number that
decomposes in~$K$ as the product~$\p_1\p_2\p_3\p_4$ where~$\p_1$,
$\p_2$ and~$\p_3$ have residue degree~$1$ and~$\p_4$ has residue
degree~$2$. We assume that the order defined in
Section~\ref{sec:primes} gives~$\p_1 < \p_2 < \p_3$. We first consider
ideals of norm~$p$. The only such ideals are~$\p_1$, $\p_2$
and~$\p_3$, and they have weight~$1$ and exponent vectors~$(1,0,0,0)$,
$(0,1,0,0)$ and~$(0,0,1,0)$. Since lexicographically we
have~$(0,0,1,0)<(0,1,0,0)<(1,0,0,0)$, the reverse lexicographic order
gives~$\p_1<\p_2<\p_3$. Now consider the ideals of norm~$p^2$, which
are $\p_4$, $\p_1^2$, $\p_2^2$, $\p_3^2$, $\p_1\p_2$, $\p_1\p_3$
and~$\p_2\p_3$. The corresponding exponent vectors are~$(0,0,0,1)$,
$(2,0,0,0)$, $(0,2,0,0)$, $(0,0,2,0)$, $(1,1,0,0)$, $(1,0,1,0)$
and~$(0,1,1,0)$, with weights~$1,2,2,2,2,2,2$. The order we defined is
therefore~$\p_4 < \p_1^2 < \p_1\p_2 < \p_1\p_3 < \p_2^2 < \p_2\p_3 <
\p_3^2$. The reader can check that the ideals of norm~$p^3$
are~$\p_4\p_1 < \p_4\p_2 < \p_4\p_3 < \p_1^3 < \p_1^2\p_2 < \p_1^2\p_3
< \p_1\p_2^2 < \p_1\p_2\p_3 < \p_1\p_3^2 < \p_2^3 < \p_2^2\p_3 <
\p_2\p_3^2 < \p_3^3$.

\subsection{Arbitrary integral ideals}

We finally define an order on the set of all nonzero integral
ideals. Again, we first order them by norm, so we only have to define
an order on nonzero ideals of the same norm. Let~$\a$ be an ideal of
norm~$N = p_1^{a_1}\dots p_k^{a_k}$ with~$p_1 < \dots <
p_k$. Then~$\a$ has a unique factorization as~$\a = \a_1\dots\a_k$
where~$\a_i$ has norm~$p_i^{a_i}$. We order ideals of norm~$N$
according to the lexicographic order of the
$k$-uple~$(\a_1,\dots,\a_k)$.

Example: let~$K$ be a number field, and assume that~$2$ decomposes
as~$\p_1\p_2\p_3$ where~$\p_1<\p_2$ have residue degree~$1$ and~$\p_3$ has
residue degree~$2$, and that~$3$ decomposes as~$\q_1\q_2$ where~$\q_1$ has
degree~$1$ and~$\q_2$ has degree~$3$. Let us order ideals of norm~$18$. Such an
ideal is uniquely a product of an ideal of norm~$2$ and an ideal of norm~$9$. The ideals
of norm~$2$ are~$\p_1 < \p_2$, and the only ideal of norm~$9$ is~$\q_1^2$. The
ideals of norm~$18$ are therefore~$\p_1\q_1^2 < \p_2\q_1^2$. Let us now order ideals
of norm~$108 = 2^2\cdot 3^3$. The ideals of norm~$4$ are~$\p_3 < \p_1^2 < \p_1\p_2 <
\p_2^2$, and the ideals of norm~$27$ are~$\q_2 < \q_1^3$. The ideals of
norm~$108$ are therefore~$\p_3\q_2 < \p_3\q_1^3 < \p_1^2\q_2 < \p_1^2\q_1^3 <
\p_1\p_2\q_2 < \p_1\p_2\q_1^3 < \p_2^2\q_2 < \p_2^2\q_1^3$.

\bibliographystyle{abbrv}

\end{document}